# New Residue Arithmetic Based Barrett Algorithms, Part I: Modular Integer Computations

## Hari Krishna Garg[1] & Hanshen Xiao[2]

**Abstract:** In this paper, we derive new computational techniques for residue number systems (RNS) based Barrett algorithm (BA). The focus of the work is an algorithm that carries out the entire computation using only modular arithmetic without conversion to large integers via the Chinese Remainder Theorem (CRT). It also avoids the computationally expensive scaling-rounding operation required in the earlier work. There are two parts to this work. First, we set up a new BA using two constants other than powers of two. Second, a RNS based BA is described. A complete mathematical framework is described including proofs of the various steps in the computations and the validity of results. Third, we present a computational algorithm for RNS based BA. Fourth, the RNS based BA is used as a basis for new RNS based algorithms for MoM and MoE. The applications we are dealing with are in the area of cryptography.

[1]Hari Krishna Garg (eleghk@nus.edu.sg), is with the Electrical & Computer Engineering Department, National University of Singapore, Singapore.
[2]Hanshen Xiao (xhs13@mails.tsinghua.edu.cn) is with the Department of Mathematics, Tsinghua University, Beijing, China.

The research work of Hari K Garg was carried out at the NUS-ZJU Sensor-Enhanced Social Media (SeSaMe) Centre. It is supported by the Singapore National Research Foundation under its International Research Centre @ Singapore Funding Initiative and administered by the Interactive Digital Media Program Office.

The work of Hanshen Xiao was supported by Tsinghua University Initiative Scientific Research Program 20141081231, the National Natural Science Foundation of China under Grant 61371104, and National High Technology Project of China under Grant 2011AA010201.

# I. Introduction

CRYPTOGRAPHY TECHNIQUES, play an important role in the security of communication and data processing systems. Instances of such cryptography techniques include **RSA** (Rivest-Shamir-Adelman), Rabin, Diffie Hellman and El Gamal. Many of these and other related techniques deal with integer arithmetic defined over a large size finite integer ring, $Z(N)$. Such an integer ring, $Z(N)$, may require several thousand bits in the binary representation of various integers in order to provide adequate level of security. A challenge in such environments is to perform the arithmetic in a computationally efficient and timely manner. The basic arithmetic operations are:

1. Modular multiplication: $C = A \cdot B \pmod{N}$, and
2. Modular exponentiation: $C = A^E \pmod{N}$.

Ordinary computations (1 and 2 as above without the mod $N$ operation) are relatively simpler to perform. However the mod $N$ operation requires division and it is a difficult operation to carry out. Further, modular exponentiation (**MoE**) is carried out via repeated use of modular multiplication (**MoM**) algorithm. Hence it is essential that an efficient algorithm be used in carrying out MoM. It is observed here that the arithmetic that may appear to be rather straightforward may turn out to be a computational bottleneck in applications such as cryptography due to involvement of large size integers. In such applications, the binary representation of integers $A$, $B$, and $N$ may turn out to be several hundred to a few thousand bits in length. The factors of $N$ may be unknown pre-empting the use of well-known Chinese Remainder Theorem (**CRT**) for the MoM and MoE computation.

Residue Number System (**RNS**) arithmetic is extensively used to express a large size computation over a finite integer ring into a number of smaller size integer rings. Computations over smaller rings can be carried out in parallel. RNS has been applied extensively in combination with computational methods, such as Barrett algorithm (**BA**) and Montgomery multiplication (**MM**), to carry out the above mentioned modulo arithmetic operations. Such methods also require base extension (**BEX**) for the RNS representation to go from one integer ring to another. In most instances, BEX turns out to be most computationally intensive part of the overall algorithm.

The major contributions of the paper are as follows. The primary and eventual objective of the work is to compute MoM and MoE as efficiently as possible with a view towards their application to cryptography. This is accomplished by first describing a new Barrett algorithm (BA) for computing the quotient associated with an integer *C* when it is divided by an integer *N*. It is assumed that *N* is a very large integer requiring several hundred to thousand bits in its binary representation. Also, its factors are unknown or at least not to be used in the BA. Second, a new RNS based Barrett algorithm (BA) is described. Third, a computationally efficient procedure for the computations underlying the RNS based BA is described. Fourth, the RNS based BA is used as a basis for new RNS based algorithms for MoM and MoE.

Mathematically, RNS represents a computation over a finite size large integer ring, $Z(M)$ into equivalent computations over a number of smaller finite size integer rings, $Z(M_i)$, $i = 1, 2, \ldots, n$. These *n* integer rings then can operate independently and hence in parallel of each other. It is desirable, therefore, to stay in the smaller size integer rings as much as possible and not carry out the computation back to the larger integer ring. This is even more urgent in cryptography when the integers can be as long as several thousand bits in their binary representation. This is also our secondary objective.

There is abundance of research materials on MM, RNS based MM, and its various extensions. A partial list of such papers is [1]-[12]. We have come across limited number of research papers on RNS based BA [13]-[14]. These two papers impose the RNS approach on the original formulation of BA meant for large integers. Hence it suffers from the major drawback of having to carry out a scaling-rounding operation twice in a single pass of the BA. The scaling-rounding operation, computing residues of $\lfloor X / 2^a \rfloor$ from residues of *X*, is computationally intensive. The algorithm described here begins with a reformulation of the basic approach to BA. It then applies RNS in a way that computations for the new BA stay within the residue arithmetic thereby avoiding scaling-rounding operation completely.

The organization of this paper is as follows. Section II provides a mathematical preliminaries on integer arithmetic, BA, MM, RNS, CRT, mixed radix systems (**MRS**), and BEX. Sections III-V are on the key contributions of the work. Mathematical structure for the new Barrett algorithm is described in Section III. The computational steps for the corresponding RNS based BA algorithm are presented in Section IV. Examples are presented to illustrate the algorithm. Section V deals with MoM and MoE algorithms that use the new Barrett algorithms. In this section, we first describe an algorithm for MoE that uses the new RNS based BA. Comparisons to existing works are also made. Section VI is on conclusions of the work.

## II. Mathematical Preliminaries

**Integer Arithmetic.** Given two integers $X$ and $A_1$, consider dividing $X$ by $A_1$ and writing

$$X = Q_1 \cdot A_1 + R_1, \tag{1}$$

where $Q_1$ is the quotient, $R_1$ is the remainder, $Q_1 \leq (X/A_1)$, and $0 \leq R_1 < A_1$. We also write (1) as

$$X \equiv R_1 \ (\text{modulo } A_1). \tag{2}$$

'Modulo' is written as 'mod' from here on. We observe that $Q_1$ and $R_1$ are unique. The process can be repeated between $Q_1$ and another integer, say $A_2$. Thus

$$Q_1 = Q_2 \cdot A_2 + R_2, \tag{3}$$

$0 \leq R_2 < A_2$. The well-known work by Euclid dating back several centuries is one of the first works dealing with such arithmetic. To take this just one step further, replace $Q_1$ in (3) by the expression in (4) to get,

$$X = Q_2 \cdot (A_2 \cdot A_1) + (R_2 \cdot A_1 + R_1). \tag{4}$$

A generalization of (4) leads to

$$X = Q_n \cdot (A_n \cdots A_1) + [R_n \cdot (A_{n-1} \cdots A_1) + \ldots + R_2 \cdot A_1 + R_1]. \tag{5}$$

Strange it may appear, (5) would be most useful in the various computations involving RNS and BEX.

**Barrett Algorithm (BA).** Given integers $A$, $B$ and $N$, $0 \leq A, B < N$ and the task of computing MoM

$$C = A \cdot B \ (\text{mod } N), \tag{6}$$

BA, consists in computing the quotient $Q$ first such that

$$X = A \cdot B \tag{7}$$

$$X = Q \cdot N + C; \ 0 \leq C < N. \tag{8}$$

Thus, $C$ is computed as

$$C = X - Q \cdot N \tag{9}$$

The computations in (8) and (9) are additional required to carry out MoM. BA is a computation of $Q$ from $X$ and $N$ by expressing $Q$ in (8) as

$$Q = \lfloor X / N \rfloor, \tag{10}$$

$\lfloor Y \rfloor$ being the floor function of $Y$. It expresses (10) as

$$Q = \lfloor X/2^a \cdot \mu / 2^{a+b} \rfloor \approx \lfloor \lfloor X/2^a \rfloor \cdot \mu / 2^b \rfloor, \qquad (11)$$

where µ is a pre-computed constant given by

$$\mu = \lfloor 2^{a+b} / N \rfloor.$$

The estimation/approximation of $Q$ in (11) is close to the actual value. The BA consists in the following three steps:

1. Compute $D = \lfloor X / 2^a \rfloor$
2. Compute $E = D \cdot \mu$
3. Compute $Q = \lfloor E / 2^b \rfloor$.

The values of scalars $a$ and $b$ can be chosen such that $Q$ as computed above is at most 2 away from the actual value of $Q$ in (10). The readers are referred to [18]-[21] for a fuller description of BA, MM, and their applications to cryptography.

**Residue number system (RNS) [2, 3].** A RNS is a finite integer ring $Z(M)$ defined by $n$ relatively co-prime moduli $M_1, M_2, \ldots, M_n$, arranged in **ascending order** without any loss of generality. The range of the RNS is given by $[0, M)$ where

$$M = \prod_{i=1}^{n} M_i. \qquad (12)$$

An integer $X \in [0, M)$ in the given RNS is represented as a length $n$ vector $\underline{X}$ via the modular computations,

$$X \leftrightarrow \underline{X} = (X_1 \, X_2 \cdots X_n), \qquad (13)$$

where

$$X_i \equiv X \,(\mathrm{mod}\, M_i), i = 1, 2, \ldots, n. \qquad (14)$$

RNS are used to express a large integer ring as a direct sum of a number of smaller integer rings. For instance, multiplication of integers $A$ and $B$ in $Z(M)$ is computed as $n$ parallel multiplications $A_i \cdot B_i \,(\mathrm{mod}\, M_i), i = 1, \ldots, n$.

**Chinese remainder theorem (CRT) [15]-[17].** Given $\underline{X}$, computation of $X$, $0 \leq X < M$, can be done via CRT, stated as

$$X \equiv \sum_{i=1}^{n} A_i \cdot X_i (\mathrm{mod}\, M_i) \cdot \left(\frac{M}{M_i}\right)(\mathrm{mod}\, M). \qquad (15)$$

The scalar $a_i$, is computed a-priori by solving the congruence

$$A_i \cdot \left(\frac{M}{M_i}\right) \equiv 1 (\bmod M_i), i = 1, 2, \ldots, n. \tag{16}$$

It is clear from (15, 16) that $\gcd(A_i, M_i) = 1, i = 1, 2, \ldots, n$. Removing mod $M$ from (15), we may also write it as

$$X \equiv \sum_{i=1}^{n} A_i \cdot X_i (\bmod M_i) \cdot \left(\frac{M}{M_i}\right) - \lambda \cdot M. \tag{17}$$

Here, the scalar $\lambda$ is an unknown and needs to be determined such that $0 \leq X < M$. The CRT computation of $X$ from $\underline{X}$ necessarily involves large integers as the dynamic range of RNS is large in cryptography applications.

**Montgomery Multiplication (MM) [1].** MM computes MoM:

$$C = A \cdot B \,(\bmod N) \tag{18}$$

$0 \leq A, B < N$, while eliminating the need to carry out 'mod $N$' operation via an intermediate computation

$$D = A \cdot B \,(\bmod R), \tag{19}$$

where $R > N$, $\gcd(N, R) = 1$, and $R \cdot R^{-1} - N \cdot N^{-1} = 1$.

The steps in MM are as follows:

1. $c \leftarrow A \cdot B$ (ordinary multiplication)
2. $d \leftarrow c \cdot N^{-1} \bmod R$
3. $e \leftarrow d \cdot N$ (ordinary multiplication)
4. $f \leftarrow c + e$
5. $g \leftarrow f / R$ (ordinary division).

The value obtained in step 5 is $(A \cdot B \cdot R^{-1}) \,(\bmod N)$ and occupies the range $[0, 2 \cdot N)$. MoE is computed via repeated application of MM where powers of $A$ and their multiplication with each other is computed via MM.

**Mixed Radix System (MRS).** MRS can also be used to compute $X$ from its residues $\underline{X}$. $X$ is represented in the following manner in an MRS:

$$X = Y_1 + Y_2 \cdot M_1 + \ldots + Y_n \cdot (M_1 \cdot M_2 \cdots M_{n-1}). \tag{20}$$

The mixed radix digits satisfy $0 \leq Y_i < M_i$. They can be computed from $\underline{X}$ using RNS to MRS conversion algorithms. For instance, see [5, 6].

**Base Extension (BEX).** Consider $\underline{X}$, residues of an integer $X$ in a RNS defined by $M$ in (8). In many situations, one requires computation of $t$ additional residues of $X$, namely $X_i \equiv X \pmod{M_i}$, $i = n + 1, \ldots, n + t$, in yet another RNS defined by $t$ additional moduli

$$M_{\mathrm{I}} = \prod_{i=n+1}^{n+t} M_i, \tag{21}$$

where $\gcd(M, M_{\mathrm{I}}) = 1$. This is called BEX. BEX tends to be expensive computationally, requiring $O(n)$ modulo operations per moduli in (21).

Broadly speaking, there are three approaches to computing BEX [7-11]. A brief description of these approaches is as follows.

**Approach 1.** As described in [7], it consists in using a redundant modulus to be used for the purposes of determining $\lambda$ in (17). Once $\lambda$ is known, (17) can then be used to compute the residues in BEX.

**Approach 2.** In this approach described in [9] and [11], one first estimates the value of $\lambda$ in (17) using the residues of $X$. Such an estimate is then used in a manner similar to approach 1.

**Approach 3.** As described in multiple places in literature on RNS including [8, 10], the residues are first be converted to mixed radix digits in MRS as per (20). The computation of BEX residues is then straightforward.

## III. A New Barrett Algorithm

A RNS based BA has been described in [13] and [14] where the authors have essentially used RNS to compute the three steps of BA, a pre-multiplication step to compute $X$ in (7), and a post-multiplication step to compute $C$ in (9). All the computations are straightforward except the first and the third steps of BA that require scaling by $2^\alpha$ in RNS. This is quite challenging in RNS and was performed via a scaling-rounding operation in [13]-[14].

Here, we first revisit the computation of $Q$ in (10). We now introduce two integers $G$ and $H$, not necessarily of the form $2^\alpha$, and approximate $Q$ in (11) as

$$Q = \lfloor X/G \cdot \mu / H \rfloor \approx \lfloor \lfloor X/G \rfloor \cdot \mu / H \rfloor, \tag{22}$$

where μ is a pre-computed constant given by

$$\mu = \lfloor G \cdot H / N \rfloor. \tag{23}$$

Since $X$ is obtained by taking the ordinary product of $A$ and $B$, both of which are assumed to be less than $N$, we have $X < N^2$. Now we derive conditions on integers $G$ and $H$ for the approximation of $Q$ in (22) to be equally close to actual $Q$ as the approximation of $Q$ in (11) is.

$\lfloor T \rfloor$ is the integral part of $T$. Consider dividing $U$ by $V$ to write $U = Q \cdot V + R$. Then we have

$$U / V = Q + R / V = \lfloor U / V \rfloor + \delta, \quad 0 \leq \delta < 1.$$

Applying the above expression to (23), we get

$$Q = \left\lfloor \frac{\left(\left\lfloor \frac{X}{G} \right\rfloor + \delta\right) \cdot \left(\left\lfloor \frac{G \cdot H}{N} \right\rfloor + \phi\right)}{H} \right\rfloor = \left\lfloor \frac{\left\lfloor \frac{X}{G} \right\rfloor \cdot \mu}{H} + \frac{\left\lfloor \frac{X}{G} \right\rfloor \cdot \phi + \left\lfloor \frac{G \cdot H}{N} \right\rfloor \cdot \delta + \delta \cdot \phi}{H} \right\rfloor, \tag{24}$$

where $0 \leq \delta, \phi < 1$. We wish the second term in the above summation to be as small as possible. To achieve that we require either (A) $\left\lfloor \frac{X}{G} \right\rfloor < H$ and $\left\lfloor \frac{G \cdot H}{N} \right\rfloor \leq H$; or (B) $\left\lfloor \frac{X}{G} \right\rfloor \leq H$ and $\left\lfloor \frac{G \cdot H}{N} \right\rfloor < H$. The difference is very subtle.

In the following, we carry out our analysis under the conditions: $0 \leq X < G \cdot H$ and $G < N < H$. Thus, $X / G < H$ and $(G \cdot H / N) < H$. We have

$$\left\lfloor \frac{X}{G} \right\rfloor \cdot \phi + \left\lfloor \frac{G \cdot H}{N} \right\rfloor \cdot \delta + \delta \cdot \phi \leq (H - 1) \cdot \phi + (H - 1) \cdot \delta + \delta \cdot \phi$$

$$= H \cdot (\phi + \delta) - (\phi + \delta) + \delta \cdot \phi.$$

Substituting in (24), we get

$$Q \leq \left\lfloor \frac{\left\lfloor \frac{X}{G} \right\rfloor \cdot \mu}{H} + \frac{H \cdot (\phi + \delta) - (\phi + \delta) + \delta \cdot \phi}{H} \right\rfloor$$

$$\leq \left\lceil \dfrac{\left\lfloor \dfrac{X}{G} \right\rfloor \cdot \mu}{H} + (\phi + \delta) \right\rceil \leq \left\lceil \dfrac{\left\lfloor \dfrac{X}{G} \right\rfloor \cdot \mu}{H} + 2 \right\rceil = \left\lfloor \dfrac{\left\lfloor \dfrac{X}{G} \right\rfloor \cdot \mu}{H} \right\rfloor + 2.$$

The above analysis and the last inequality leads to the conclusion summarized in the following theorem:

**Theorem 1.** If the integers $G$ and $H$, $G < H$, satisfy the conditions

$$N^2 \leq G \cdot H, \text{ and } G < N, \tag{25}$$

then the estimate/approximate value of quotient is at most 2 away from the actual quotient.

Thus the resulting remainder computed in (9) satisfies,

$$0 \leq C < 3 \cdot N. \tag{26}$$

A final correction may be required to reduce $C$ to the range $[0, N)$. It is interesting to note that $\gcd(G, H) = 1$ is not required. All the above analysis results in the new BA summarized in the following.

**A New Barrett Algorithm for Computing $A \cdot B \bmod N$**

Input: $A, B, N, G, H; 0 \leq A, B < N, G < N, N^2 \leq G \cdot H$.

Output: $Q$ (an estimate of quotient when $A \cdot B$ is divided by $N$)

**Step 1.** Compute $\mu = \lfloor G \cdot H / N \rfloor$ (one-time)

**Step 2.** Compute $X = A \cdot B$ (ordinary multiplication)

**Step 3.** Compute $D = \lfloor X / G \rfloor$ (quotient)

**Step 4.** Compute $E = D \cdot \mu$ (ordinary multiplication)

**Step 5.** Compute $Q = \lfloor E / H \rfloor$ (quotient)

Once $Q$ is computed, the algorithm may proceed further to compute the remainder $X \bmod N$ as

**Step 6.** $C = X - Q \cdot N$ (remainder, ordinary multiplication)

As observed earlier, $0 \leq C < 3 \cdot N$. In many cases, a further step may be required to reduce $C$ further to a range $[0, N)$. In many others, such a step may not be needed.

The conditions that are required for the integers $G$ and $H$ are general ($G < N, N^2 \leq G \cdot H$) and leave door open to a wide range of possibilities. The corresponding computational steps can also be vastly different. For instance, $G = 1$ and $N^2 \leq H$ is valid. In such a case, Step 2 as described above is not required as $D = X$ when $G = 1$. Other sets of conditions may also be derived from (24). For instance, $G = N$ and $N \leq H$ can also be used. The analysis is straightforward to carry out.

It is also possible to use the above analysis to explore other possibilities. For instance, we may wish to have $0 \leq C < 2 \cdot N$. In such as case, we may set $G < 0.5 \cdot N$ and $2 \cdot N^2 \leq G \cdot H$. These values lead to a value of $Q$ that is at most 1 away from its actual value. For now, we make these observations to make the readers aware of the general nature of the analysis and the various possibilities for the new BA. Also, we are quite sure that $G$ and $H$ will never assume the same value as $N$. Otherwise the original computation $A \cdot B$ mod $N$ can be performed instead.

**Example 1.** Let $N = 21$, $A = 20$, $B = 19$. Thus, $C = X$ mod $21 = 380$ mod $21 = 2$. The conditions to be satisfied are $G < 21$ and $441 < G \cdot H$. Choosing $G = 20$, $H = 24$, we have, $\mu = 22$, $X = 380$, $D = 19$, $E = 418$, $Q = 17$, and $C = 23$. Since $C > 21$, we need to subtract 21 from this value to get the final result $C = 2$.

**Example 2.** We work with the same values for $N$, $A$, and $B$ as in example 1. Choosing $G = 1$, $H = 600$, we have, $\mu = 28$, $X = 380$, $D = 380$, $E = 10{,}640$, $Q = 17$, and $C = 23$.

**Example 3.** We work with the same values for $N$, $A$, and $B$ as in example 1. Choosing $G = 10$, $H = 89$, we have, $\mu = 42$, $X = 380$, $D = 38$, $E = 1{,}596$, $Q = 17$, and $C = 23$.

As described here, BA computes $C = A \cdot B$ mod $N$, where $0 \leq A, B < N$ and $0 \leq C < 3 \cdot N$. However, in many situations, for instance when we use MoM recursively to compute MoE, we may wish to use BA repeatedly. In those cases, $C$ will become input to MoM in the next step. Hence it may be desirable to define BA in a manner that the inputs $A$ and $B$ and the output $C$ occupy the same range $[0, 3 \cdot N)$. This can be achieved by replacing the first condition $N^2 \leq G \cdot H$ by $9 \cdot N^2 \leq G \cdot H$ in the algorithm described above. In a similar manner, if we wish to compute $C = A \cdot B$ mod $N$, where $0 \leq A, B, C < 2 \cdot N$, then we may select $G$ and $H$ such that $G < 0.5 \cdot N$ and $8 \cdot N^2 < G \cdot H$. This results in a slight increase in the range of values admissible for $G$ and $H$.

Table 1 shows a summary of the discussion thus far.

| **Case 1:** $0 \leq A, B < N$ | **Case 2:** $0 \leq A, B < 3 \cdot N$ | **Case 3:** $0 \leq A, B < N$ | **Case 4:** $0 \leq A, B < 2 \cdot N$ |
|---|---|---|---|
| $0 \leq C < 3 \cdot N$ | $0 \leq C < 3 \cdot N$ | $0 \leq C < 2 \cdot N$ | $0 \leq C < 2 \cdot N$ |
| $G < N, N^2 < G \cdot H$ | $G < N, 9 \cdot N^2 < G \cdot H$ | $G < 0.5 \cdot N, 2 \cdot N^2 < G \cdot H$ | $G < 0.5 \cdot N, 8 \cdot N^2 < G \cdot H$ |
| $G_{max} \approx N, H_{min} \approx N$ | $G_{max} = N, H_{min} \approx 9 \cdot N$ | $G_{max} \approx 0.5 \cdot N, H_{min} \approx 4 \cdot N$ | $G_{max} \approx 0.5 \cdot N, H_{min} \approx 16 \cdot N$ |

**Table 1. Conditions on $G$ and $H$ for $C = A \cdot B \bmod N$.**

In the next section, we turn to the various computations for BA algorithm described above in RNS. We note that there is no scaling-rounding operation to be carried out in the RNS arithmetic based new BA described here.

### IV. RNS Based New Barrett Algorithm

We now turn to a RNS based BA for the computation in Steps 1-6. All the integers need to be expressed in residue form over multiple moduli that constitute a suitable finite integer ring $Z(M)$. Step 1 involves a one-time computation, so that is straightforward to map to a RNS. Steps 2, 4, and 6 involve ordinary multiplication of two integers and that is also straightforward to map to a RNS. Steps 3 and 5 require computation of quotient, or alternately stated scaling, by integers $G$ and $H$ respectively. As we will see in the following, in order for us to compute these two steps solely with residue arithmetic, both $G$ and $H$ must be a factor of the moduli that constitute the RNS used. Also, since the BA is expected to be used recursively for carrying out MoE, we would like to use the same RNS in all Steps 1-6. In addition, for RNS residues to correspond to the actual integer values in Steps 1-6, the product of the moduli of the RNS must exceed the maximum value of the integer possible at each step.

Starting with $0 \leq A, B \leq N$, the maximum integer value possible in each step is as follows:

Step 1. $\mu \leq G \cdot H / N$ (Integer part if not an integer)

Step 2. $X \leq (N-1)^2$

Step 3. $D \leq (N-1)^2 / G$ (Integer part if not an integer)

Step 4. $E \leq H \cdot (N-1)^2 / N < H \cdot N$

Step 5. $Q \leq N - 1$

Step 6. $Q \cdot N \leq (N - 1) \cdot N.$

Based on the above bounds on integer values and the conditions in Theorem 1, we have the following conditions for the RNS based on $Z(M)$ to be used the RNS BA:

1. $N^2 \leq G \cdot H$

2. $G \leq N$

3. $G \,|\, M$

4. $H \,|\, M$

5. $H \cdot N < M.$

For the case described earlier when estimated $Q$ is at most 1 away from the actual value, we have the following conditions for the RNS based on $Z(M)$ to be used the RNS BA:

1. $2 \cdot N^2 \leq G \cdot H$

2. $G < 0.5 \cdot N$

3. $G \,|\, M$

4. $H \,|\, M$

5. $H \cdot N < M.$

It is interesting to note that $\gcd(G, H) = 1$ is not required. Further, even in all cases including when $0 \leq A, B, C < 3 \cdot N$ and $0 \leq A, B, C < 2 \cdot N$, the largest integer that needs to be computed is $E$ in step 4. Given that $E = D \cdot \mu = \lfloor X / G \rfloor \cdot \lfloor G \cdot H / N \rfloor \leq \lfloor (X / G) \cdot (G \cdot H / N) \rfloor = \lfloor X \cdot H / N \rfloor < X \cdot H / N$. Thus, the condition that

$$X \cdot H / N < M \qquad (27)$$

must be valid in all instances.

A large number of possibilities become apparent now. We will go about selecting $G$ and $H$ first that satisfy the conditions. Then the integer $M$ is constructed such that $\text{lcm}(G, H) \,|\, M$. Finally, if needed further suitable residues are included in $M$ to satisfy the last condition $H \cdot N < M$.

The computational procedure for steps 3 and 5 require the computation of quotients in the RNS. Such an algorithm is described in Appendix A. Since it requires modulo inverses, it can work only when the various moduli are co-prime. Also, the result of $\lfloor X / G \rfloor$ would only be known in terms of residues corresponding to the moduli that constitute $M / G$. Thus, for the quotient to be computed, we require $\gcd(G, M / G) = 1$. Similarly, it is also required that $\gcd(H, M / H) = 1$. Again, $\gcd(G, H) = 1$ is not required.

The above description of the computation of quotient residues also brings out another interesting aspect. The computations in Step 1 and 2 are carried out over $Z(M)$. After Step 3 computation is done, the quotient $\lfloor X / G \rfloor$ is available in $Z(M / G)$. Hence we need to carry out BEX to expand the quotient residues back to $Z(M)$. Similarly, we need to carry out another BEX to expand the quotient residues of $\lfloor E / H \rfloor$ computed in $Z(M / H)$ to $Z(M)$. Such a BEX algorithm is described in Appendix B [22].

**The new RNS based Barrett Algorithm**

**Given**

RNS defined over $Z(M)$, residues of $N, G, H$ in $Z(M)$

In step 2a, the first $a$ factors of $M$ correspond to $G$ while in step 4a, the first $b$ correspond to $H$. This is assumed without any loss in generality.

**Input**

Residues of $A$ and $B$, that is, $(A_i, B_i) \equiv (A, B) \bmod M_i$, $i = 1, \ldots, n$.

**Pre-computational Step**

**Step 1.** Compute $\mu_i$, $i = 1, \ldots, n$, residues of $\mu = \lfloor G \cdot H / N \rfloor$ in $Z(M)$.

**Computation Steps**

**Step 2. Modulo multiplication.** Compute $X_i \equiv A_i \cdot B_i \pmod{M_i}$, $i = 1, \ldots, n$.

**Step 3a. Quotient computation.** Let $M_I = G$ and $M_{II} = M / M_I$. Compute residues $D_i$, $i = a + 1, \ldots, n$, from residues $X_i$, $i = 1, \ldots, n$, $G_i = M_i$, $i = 1, \ldots, a$.

**Step 3b. BEX.** Use BEX on residues $D_i$, $i = a + 1, \ldots, n$, to get $u$ residues $D_i$, $i = 1, 2, \ldots, a$.

**Step 4. Modulo multiplication.** Compute $E_i \equiv D_i \cdot \mu_i$, $i = 1, \ldots, n$.

**Step 5a. Quotient computation.** Let $M_I = H$ and $M_{II} = M / M_I$. Compute residues $Q_i$, $i = b + 1, \ldots, n$, from residues

$E_i$, $i = 1, \ldots, n$, $H_i = M_i$, $i = 1, \ldots, b$.

**Step 5b. BEX.** Use BEX on residues $Q_i$, $i = b + 1, \ldots, n$, to get $b$ residues $Q_i$, $i = 1, 2, \ldots, b$.

**Step 6. Remainder computation.** Compute $C_i \equiv X_i - Q_i \cdot N_i$, $i = 1, \ldots, n$.

**END**

**Example 4.** Let $N = 21$, $G = 20$, $H = 28$. We have the conditions $20 \mid M$, $28 \mid M$, and $588 < M$. We select $M = 1{,}540 = 4 \cdot 5 \cdot 7 \cdot 11$, $\mathbf{m} = (4\ 5\ 7\ 11)$, for the computation. The conditions $\mathrm{lcm}(G, H) \mid M$, $\gcd(G, M/G) = \gcd(H, M/H) = 1$ are satisfied. Let $A = 20$, $B = 19$. The answer is $20 \cdot 19 \bmod 21 = 2$. The various residue vectors that the RNS based new BA processes are as follows.

Input $\underline{\mathbf{A}} = (0\ 0\ 6\ 9)$, $\underline{\mathbf{B}} = (3\ 4\ 5\ 8)$

Step 1. $\underline{\boldsymbol{\mu}} = (2\ 1\ 5\ 4)$

Step 2. $\underline{\mathbf{X}} = (0\ 0\ 2\ 6)$

Step 3a. $M_a = 20 = m_1 \cdot m_2$, $\underline{\mathbf{D}} = (*\ *\ 5\ 8)$.

Step 3b. $\underline{\mathbf{D}} = (3\ 4\ 5\ 8)$

Step 4. $\underline{\mathbf{E}} = (2\ 5\ 4\ 10)$

Step 5a. $M_a = 28 = m_1 \cdot m_3$, $\underline{\mathbf{Q}} = (*\ 2\ *\ 6)$

Step 5b. $\underline{\mathbf{Q}} = (1\ 2\ 3\ 6)$

Step 6. $\underline{\mathbf{C}} = (3\ 3\ 2\ 1)$.

The algorithm calculates 23 as the answer which is $N$ away from the correct answer of 2.

The above described RNS based BA computes $C = A \cdot B \bmod N$, where $0 \leq A, B < N$ and $0 \leq C < 3 \cdot N$. However, in many situations we may wish to use RNS BA repeatedly, for instance in computing modular exponentiation. Hence it may be desirable to define RNS BA in a manner that the inputs $A$ and $B$ and the output $C$ occupy the same range $[0, 3 \cdot N)$. This can be achieved by replacing the first condition $N^2 \leq G \cdot H$ by $9 \cdot N^2 \leq G \cdot H$ and the last condition $H \cdot N < M$ by $9 \cdot H \cdot N < M$ in the algorithm described above. Further, we may also wish to limit $A$, $B$, and $C$ to the range $[0, < 2 \cdot N)$. This may be achieved by replacing the first condition $2 \cdot N^2 \leq G \cdot H$ by $8 \cdot N^2 \leq G \cdot H$ and the last condition $H \cdot N < M$ by $4 \cdot H \cdot N < M$. We note that this results in a slight increase in the size of the RNS. However, in cryptography, the numbers are large and this increase would be infinitesimally small.

It is also possible to simplify the BEX computations in steps 3b and 5b by using one redundant residue in $M$, say $M_{n+1}$ such that $M_{n+1} \geq \max(n - a, n - b)$. In this case, the well-known method in [7], also described as Approach 1 in Section II, can be used in a straightforward manner. It will lead to considerable computational savings for BEX in steps 3b and 5b at the expense of very slight increase in complexity of other steps.

We end this section that there is a wide variety of algorithms can be obtained by simply examining various possibilities for the values of $G$ and $H$. The range associated with the resulting quotient can also be broadened or narrowed. It is indeed possible to select $G$ and $H$ such that $G \mid H$. This may lead to simplification in the hardware implementation as the computation of $\lfloor X / G \rfloor$ in step 3 and $\lfloor E / H \rfloor$ in step 5 may be done using very similar architecture and parameters.

## V. Further Analysis

In this section, we first describe an algorithm for MoE that uses the new RNS based BA. In this case, we are assuming that the input integers $A$ and $B$ and the output integer $C$ occupy the same range, say $[0, 3 \cdot C)$. We term such an algorithm Barrett modular multiplication (**BMM**) algorithm. Pseudocode for BMM and MoE algorithms are shown in the following.

**BMM Algorithm (A, B)**

**Input:** Residue vectors $\underline{A}$, $\underline{B}$ in $Z(M)$, $0 \leq A, B < 3 \cdot N$.

**Output:** Residue vector $\underline{C}$ in $Z(M)$ such that $C = A \cdot B \pmod{N}$, $0 \leq C < 3 \cdot N$.

Step 1. Precompute residue vectors $\underline{\mu}$ and $\underline{N}$ in $Z(M)$.

Step 2. Compute $\underline{X} = \underline{A} \cdot \underline{B}$ in $Z(M)$.

Step 3. Compute residues $\underline{D}$ for quotient $D = \lfloor X / G \rfloor$ in $Z(M / G)$.

Step 3a. Compute BEX to get residues $\underline{D}$ in $Z(M)$.

Step 4. Compute $\underline{E} = \underline{D} \cdot \underline{\mu}$ in $Z(M)$.

Step 5. Compute residues **Q** for quotient $Q = \lfloor E/H \rfloor$ in $Z(M/H)$.

Step 5a. Compute BEX to get residues **Q** in $Z(M)$.

Step 6. Compute $\underline{C} = \underline{X} - \underline{Q} \cdot \underline{N}$ in $Z(M)$.

**BMM based MoE**

**Input:** Residue vector $\underline{X}$ in $Z(M)$ and $E$, $E = \sum_{i=0}^{k} e_i \cdot 2^i$, $0 \leq X < 3 \cdot N$.

**Output:** Residues $\underline{Y}$ for $Y = X^E \pmod{N}$ in $Z(M)$, $0 \leq Y < 3 \cdot N$

1. If $e_0 = 1$

   $\underline{Y} \leftarrow \underline{X}$

   else

   $\underline{Y} \leftarrow \underline{1}$

2. For $j = 1$ to $k$ do

   $\underline{X} \leftarrow \text{BMM}(\underline{X}, \underline{X})$

   If $e_j = 1$ then

   $\underline{Y} \leftarrow \text{BMM}(\underline{Y}, \underline{X})$

   end If

   end For.

Here, $\underline{1}$ denotes the vector of all 1s. It is straightforward to describe BMM and BMM based MoE when $0 \leq A, B, C < 2 \cdot N$. Comapred to the MoE that uses existing RNS based MM [7, 8, 9, 10], the proposed Barrett based MoE directly returns $X^E \bmod N$ rather than $(X^E \cdot R^{-1}) \bmod N$. Further, there is greater flexibility in selecting $G$ and $H$ in RNS based Barrett algorithm.

## VI. Conclusions

In this work, a new Barrett algorithm is described for computing modulo multiplication $A \cdot B \pmod{N}$ and modulo exponentiation $A^E \bmod N$. A residue number system based version of the new Barrett algorithm is also described that uses only residue arithmetic thereby avoiding the scaling-rounding operation that is computationally intense. Both algorithms as described here are a first from our understanding of the research literature. The previously known Barrett algorithms use powers to 2 to scale the various computations. Similarly, the previously known residue number system based Barrett algorithm requires expensive scaling-rounding operations. We also described the application of the new residue number system based Barrett algorithm to computing modular exponentiation. Comparisons with existing works are also made. We have also completed our research work on suitable architectures with application in cryptography as our main focus.

**References**

[1] PL Montgomery, "Modular Multiplication Without Trial Division," *Mathematics of Computation*, vol 44, no 170, pp 519-521, 1985.

[2] PV Ananda Mohan, Residue Number Systems, Algorithms and Architectures, *Springer International Series in Engineering & Computer Science*, 2002.

[3] A Omondi & B Premkumar, Residue Number Systems, Theory & Implementation, *Imperial College Press, Advances in Computer Science and Engineering*, 2007.

[4] J-C Bajard & T Plantard, "RNS Bases and Conversions," *Proceedings SPIE – International Society for Optical Engineering,* vol 5559, pp 60-69, 2004.

[5] AP Shenoy & R Kumaresan, "Fast Base Extension using a Redundant modulus in RNS," *IEEE Trans. on Computers,* vol 38, pp 292-297, 1989.

[6] J-C Bajard, L-S Didier & P Kornerup, "Modular Multiplication and Base Extensions in Residue Number Systems," *IEEE Symposium on Computer Arithmetic,* pp 59-65, 2001.

[7] S Kawamura, M Koike, F Sano & A Shimbo, "Cox-Rower Architecture for Fast Parallel Montgomery Multiplication," Advances in Cryptography – EUROCRYPT 2000, Springer Lecture Notes in Computer Science, vol 1807, pp 523-538, 2000.

[8] F Gandino, F Lamberti, G Paravati, J-C Bajard & P Montuschi, "An Algorithmic and Architectural Study on Montgomery Exponentiation in RNS," *IEEE Trans. on Computers,* vol 61, pp 1071-1083, 2012.

[9] KC Posch & R Posch, "Modulo Reduction in Residue Number Systems," *IEEE Transactions on Parallel and Distributed Systems*, vol 6, pp 449-454, 1995.

[10] J Bajard, L-S Didier, & P Kornerup, "An RNS Montgomery modular multiplication algorithm," *IEEE Transactions on Computers*, vol 47, pp 766-776, 1998.

[11] H Nozaki, M Motoyama, A Shimbo, & S Kawamura, "Implementation of RSA Algorithm Based on RNS Montgomery Multiplication," *Cryptographic Hardware and Embedded Systems - CHES 2001*, vol 2162, pp 364-376, 2001.

[12] KA Gbolagade & SD Cotofana, "An O(*n*) Residue Number System to Mixed Radix Conversion technique," *IEEE International Symposium on Circuits and Systems, ISCAS*, pp 521-524, 2009.

[13] D Schinianakis & T Stouraitis, "An RNS Barrett Modular Multiplication Architecture," *IEEE International Symposium on Circuits & Systems, ISCAS*, pp 2229-2232, 2014.

[14] D Schinianakis & T Stouraitis, "An RNS Modular Multiplication Architecture," *IEEE International Conference on Electronics, Circuits & Systems (ICECS)*, pp 958-961, 2013.

[15] C Ding, D Pei, & A Salomaa, Chinese Remainder Theorem, Applications in Computing, Coding, and Cryptography, *World Scientific*, 1996.

[16] H Krishna, B Krishna, K-Y Lin & J-D Sun, Computational number theory and digital signal processing: Fast algorithms and error control techniques, *CRC Press*, Boca Raton, USA, 1994.

[17] NS Szabo & RI Tanaka, Residue arithmetic and its applications to computer technology, *McGraw-Hill*, New York, 1967.

[18] A Menezes, P van Oorschot, & S Vanstone, Handbook of Applied Cryptography, *CRC Press*, 1996.

[19] LC Washington, Elliptic Curves: Number Theory & Cryptography, *CRC Press*, 2008.

[20] D Hankerson, A Menezes & S Vanstone, Guide to Elliptic Curve Cryptography, *Springer-Verlag New York*, 2004.

[21] Z Cao, R Wei, & X Lin, "A fast modular reduction method," *IACR Cryptology ePrint Archive*, 2014.

[22] HK Garg and H Xiao, "A New BEX Technique in Residue Systems," submitted to *IEEE Signal Processing Letters*, Jan 2016.

# APPENDIX A: Computing Quotient in RNS

**Problem Formulation.** Given residues $x_i$ of an integer $X$ in $Z(M)$, $X_i \equiv X \pmod{M_i}$, $i = 1, \ldots, n$, $M = M_I \cdot M_{II}$,

$$M_I = \prod_{i=1}^{a} M_i,$$

$$M_{II} = \prod_{i=a+1}^{n} M_i.$$

Compute the residues of the quotient $Q$ when $X$ is divided by $M_I$, $0 \leq Q < M_{II}$.

We now revisit integer arithmetic introduced in Section II and use it to describe the algorithm for computing the residues of $Q$. Consider (1) when $X$, $A_1$ and $R_1$ are known. We may compute $Q_1$ as

$$Q_1 = (X - R_1) \cdot A_1^{-1} \tag{A1}$$

Again, when $Q_1$, $A_2$ and $R_2$ are known we may compute $Q_2$ as

$$Q_2 = (Q_1 - R_2) \cdot A_2^{-1} \tag{A2}$$

This can be carried out recursively to finally compute $Q_a$ as

$$Q_a = (Q_{a-1} - R_a) \cdot A_a^{-1} \tag{A3}$$

The representation of $X$ in (5) is still valid and is reproduced below for completeness.

$$X = Q_a \cdot (A_a \cdots A_1) + [R_a \cdot (A_{a-1} \cdots A_1) + \ldots + R_2 \cdot A_1 + R_1]. \tag{A4}$$

Let us apply the above integer arithmetic in (A1)-(A4) to the RNS defined over the integer ring $Z(M)$. We are given moduli $M_i$ and corresponding residues $X_i$, $i = 1, \ldots, n$. We now set

$$A_i = M_i. \tag{A5}$$

Thus,

$$R_1 = X_1. \tag{A6}$$

This leads to,

$$Q_1 = (X - X_1) \cdot M_1^{-1} \tag{A7}$$

Since $X$ is expressed in terms of its residues and $M_1^{-1}$ exists only mod $M_i$, $i = 2, \ldots, n$, we compute residues of $Q_1$ in (A7) by taking mod $M_i$, $i = 2, \ldots, n$, of both sides. This results in

$$Q_{1,i} \equiv (X_i - X_1) \cdot M_1^{-1} \pmod{M_i}, \quad i = 2, \ldots, n. \tag{A8}$$

Since, $Q_1 \leq X / M_1$, it is uniquely expressed in terms of its residues $Q_{1,i}$, $i = 2, \ldots, n$. After the first iteration in (A1), it is seen that

$$R_2 \equiv Q_2 \pmod{M_2} = Q_{1,2}. \tag{A9}$$

Again, expressing (A9) in residue form, we compute residues of $Q_2$ by taking mod $M_i$, $i = 3, \ldots, n$, of both sides. This results in

$$Q_{2,i} \equiv (Q_{1,i} - Q_{1,2}) \cdot M_2^{-1} \pmod{M_i}, i = 3, \ldots, n. \tag{A10}$$

Again $Q_2 \leq Q_1 / M_2$ and hence uniquely expressed in terms of its residues $Q_{2,i}$, $i = 3, \ldots, n$.

This process is carried out for $a$ iterations computing residues of $Q_k$ or $Q_{k,i}$, $i = k + 1, \ldots, n$, at the $k$-th iteration, $k = 1, \ldots, a$. Up on conclusion of the $a$ iterations, we have (A4) as

$$X = Q_a \cdot (M_a \cdots M_1) + [R_a \cdot (M_{a-1} \cdots M_1) + \ldots + R_2 \cdot M_1 + R_1]. \tag{A11}$$

Thus, $Q_a$ is the remainder obtained by dividing $X$ by $M_I$ expressed in terms of its residues $Q_{a,i}$, $i = a + 1, \ldots, n$.

## APPENDIX B: Computing Base Extension in RNS

**Problem Formulation.** BEX computational problem is:

Given

1. two RNS defined over $Z(M_I)$ and $Z(M_{II})$,
2. Residues of $X$ in $Z(M_I)$, that is, $X_i \equiv X \pmod{M_i}$, $i = 1, \ldots, a$.

Compute

$n - a$ residues of $X$ over $Z(M_{II})$, that is, $X_i \equiv X \pmod{M_i}$, $i = a + 1, \ldots, n$.

We wish not to involve CRT in BEX computation either via (15) or (17).

The representation of $X$ in (5) is still valid and is reproduced below for completeness.

$$X = Q_n \cdot (A_n \cdots A_1) + [R_n \cdot (A_{n-1} \cdots A_1) + \ldots + R_2 \cdot A_1 + R_1]. \tag{B1}$$

Let us apply the above integer arithmetic to the RNS defined over the integer ring $Z(M)$. We are given moduli $M_i$ and corresponding residues $X_i$, $i = 1, \ldots, a$. We now extend the moduli set to include moduli $M_i$, $i = a + 1, \ldots, n$. Since the corresponding residues are unknown (and need to be computed), we extend the residue set to include residues $X_i$, $i = a + 1, \ldots , n$, by choosing any arbitrary values (can be all 0s) for them. In the current formulation, we set

$$A_i = M_i. \tag{B2}$$

Thus,

$$R_1 = X_1. \tag{B3}$$

This leads to,

$$Q_1 = (X - X_1) \cdot M_1^{-1} \tag{B4}$$

Since $X$ is expressed in terms of its residues and $M_1^{-1}$ exists only mod $M_i$, $i = 2, \ldots , n$, we compute residues of $Q_1$ in (B4) by taking mod $M_i$, $i = 2, \ldots , n$, of both sides. This results in

$$Q_{1,i} \equiv (X_i - X_1) \cdot M_1^{-1} \pmod{M_i}, i = 2, \ldots , n. \tag{B5}$$

Since, $Q_1 \leq X / M_1$, it is uniquely expressed in terms of its residues $Q_{1,i}$, $i = 2, \ldots , n$. After the first iteration in (B5), it is seen that

$$R_2 \equiv Q_2 \pmod{M_2} = Q_{1,2}. \tag{B6}$$

Again, expressing (B6) in residue form, we compute residues of $Q_2$ by taking mod $M_i$, $i = 3, \ldots , n$, of both sides. This results in

$$Q_{2,i} \equiv (Q_{1,i} - Q_{1,2}) \cdot M_2^{-1} \pmod{M_i}, i = 3, \ldots , n. \tag{B7}$$

Again $Q_2 \leq Q_1 / M_2$ and hence uniquely expressed in terms of its residues $Q_{2,i}$, $i = 3, \ldots , n$.

This process is carried out for $a$ iterations computing residues of $Q_k$ or $Q_{k,i}$, $i = k + 1, \ldots , n$ at the $k$-th iteration, $k = 1, \ldots , a$. Up on conclusion of $a$ iterations, we have

$$X = Q_a \cdot (M_a \cdots M_1) + [R_a \cdot (M_{a-1} \cdots M_1) + \ldots + R_2 \cdot M_1 + R_1]. \tag{B8}$$

A quick observation of (B8) reveals the term in [.] is the mixed-radix digits of $X$ corresponding to residues $X_i$, $i = 1, \ldots , a$. Also, residues of $Q_a$ are known for $M_i$, $i = a + 1, \ldots , n$. Thus we have

$$X \text{ in given MRS} = X - Q_a \cdot (M_a \cdots M_1). \tag{B9}$$

Expressing (B9) in residue form, we get residues in RNS over $Z(M_{II})$ as

$$X_i \equiv (X \text{ in given MRS}) \bmod M_i \equiv X_i - Q_{a,i} \cdot M_1 \pmod{M_i}, i = a + 1, \ldots , n. \tag{B10}$$

Thus, residues up on BEX are computed as per (B10). The Steps of the BEX algorithm are as follows:

**The new BEX algorithm**

**Input:**

Two RNS defined over $Z(M_I)$ and $Z(M_{II})$;

Residues of $X$ over $Z(M_I)$, that is $X_i \equiv X \pmod{M_i}$, $i = 1, \ldots, a$.

**Initialization:**

Assign arbitrary values (can be all 0s) to $X_i$, $i = a + 1, \ldots, n$; and initialize $Q_0$ as

$$Q_{0,i} = X_i, i = 1, \ldots, n. \tag{B11}$$

**Computational Steps:** For $k = 1, \ldots, a$, compute

$$Q_{k,i} \equiv (Q_{k-1,i} - Q_{k-1,k}) \cdot M_k^{-1} \pmod{M_i}, i = k + 1, \ldots, n. \tag{B12}$$

**Final Computational Step:** For $i = a + 1, \ldots, n$, compute the BEX residues as

$$X_i \equiv X \pmod{M_i} \equiv X_i - Q_{a,i} \cdot M_a \pmod{M_i}, i = a + 1, \ldots, n. \tag{B13}$$

END

We assume that all modulo products and inverses associated with various moduli are precomputed and stored. It will be clear to the readers that all the quantities computed by the new BEX algorithm are unique and initial arbitrary values assigned to $X_i$, $i = a + 1, \ldots, n$, have no bearing on the BEX values. This is stated as a theorem in the following.

**Theorem.** The residue values computed via the BEX algorithm described here leads to the correct residues for all the arbitrarily assigned initial values to $X_i$, $i = a + 1, \ldots, n$.

Due to their wider appeal, the materials of this appendix have also been submitted separately for publication [22]. We end this appendix by stating that one may also use any of the algorithms for BEX as described in Section II.